\documentclass[a4paper,12p]{article}
\title{Equiaffine structure and conjugate Ricci-symmetry of a statistical manifold}
\author{\textrm{Chol Rim Min , Won Hak Ri , Hyong Chol O} $~~$ $~~$\\  \\
                  {\small\textit{Faculty of Mathematics }}\\
                        {\small\textit{Kim Il Sung University, Pyongyang, D. P. R. Korea}}\\   
}
    
\date{January, 2013}  
\usepackage{latexsym}   
\usepackage[pdftex]{graphicx}
\usepackage{amsmath}
\newtheorem{Thot}{Theorem}[section]
\newtheorem{Prop}{Proposition}[section]
\newtheorem{Coro}{Corollary}[section]
\newtheorem{Def}{Definition}[section]
\newtheorem{lemm}{Lemma}[section]
    
\begin{document}
\maketitle      
\begin{abstract}   
A condition for a statistical manifold to have an equiaffine structure is studied. The facts that dual flatness and conjugate symmetry of a statistical manifold are sufficient conditions for a statistical manifold to have an equiaffine structure were  obtained in [2] and [3]. In this paper, a fact that a statistical manifold, which is conjugate Ricci-symmetric, has an equiaffine structure is given, where conjugate Ricci-symmetry is weaker condition than conjugate symmetry.
A condition for conjugate symmetry and conjugate Ricci-symmetry to coincide is also given.
\\


$2010 Mathematics Subject Classiﬁcation$. Primary 53A15; Secondary 53B05, 53C44 \\

\end{abstract}
\section{Introduction}
\quad Information geometry proposed by Amari has been applied to various fields of mathematical sciences such as mathematical statistics, statistical physics, neural networks and information theory. Nowadays, the relation between Information geometry and Bayesian statistics has been studied. 
 It is known that Jeffreys prior in Bayesian statistics is the volume element of the Fisher metric on the manifold of probability distributions.\\
 Let $M$ be an $n$-dimensional manifold, and $\nabla$  be a torsion-free affine connection on $M$. Let $\omega$ be a volume form on $M$, that is, $\omega$ is an $n$-form which does not vanish everywhere on $M$.\\
If $\nabla\omega=0$ , we say that $\nabla$ is a (local) equiaffine connection, and, $\omega$ is a parallel volume form on $M$. In this case $(\nabla,\omega)$ is called a (local) equiaffine structure on $M$.\\
 The parallel volume form with respected to a given connection $\nabla$ is determined uniquely up to constant multiplications. For more details, see [5].\\
From a viewpoint of differential geometry, Jeffreys prior in Bayesian statistics is a parallel volume form with respected to $\stackrel{\circ}\nabla$, where $\stackrel{\circ}\nabla$ is a Levi-Civita connection of Fisher metric.\\
For all $\alpha\in\mathbf{R}$, a connection $\nabla^{(\alpha)}$ is defined by
\begin{equation*}
	\nabla^{(\alpha)}=\frac{1+\alpha}{2}\nabla+\frac{1-\alpha}{2}\stackrel{*}\nabla
\end{equation*}
where $\nabla$  and $\stackrel{*}\nabla$  are dual connections on $M$. For more details, see [1].
Jeffreys prior was generalize to $\alpha$-parallel prior with respected to $\alpha$- connection for all $\alpha\in\bf R$ in [3] and the researches on them have been done in [2], [3] and [4].
On the other hand, it is known that a necessary and sufficient condition for $\nabla$ to be equiaffine is that Ricci curvature tensor of $\nabla$ is symmetric. For more details, see [5].\\
Some conditions for $\nabla^{(\alpha)}$ to be equiaffine have been studied for all $\alpha\in\mathbf{R}$.\\
It was given that dual flatness of a statistical manifold $(M, g, \nabla, \stackrel{*}\nabla)$ is a sufficient condition for $\nabla^{(\alpha)}$ to be equiaffine in [4], and that conjugate symmetry of a statistical manifold $(M, g, \nabla, \stackrel{*}\nabla)$ is also so in [2]. For more details of conjugate symmetry of a statistical manifold, see [6].\\
This paper gives a fact that conjugate Ricci symmetry of a statistical manifold $(M, g, \nabla)$ is a sufficient condition for a statistical manifold $(M, g, \nabla)$ to have an equiaffine structure, where conjugate Ricci-symmetry of a statistical manifold is weaker than conjugate symmetry of a statistical manifold.\\
In section 2, the properties of curvature tensor field of a statistical manifold $(M, g, \nabla)$  are given, which is also preliminary for defining of conjugate Ricci-symmetry of a statistical manifold $(M, g, \nabla)$  in section 3.\\
In section 3, conjugate Ricci-symmetry of a statistical manifold $(M, g, \nabla)$  is defined, and some properties associated with an equiaffine structure of a statistical manifold are given.\\
In section 4, a statistical manifold with a recurrent metric is defined, and it is shown that if some conditions are satisfied, conjugate Ricci-symmetry of it coincides with conjugate symmetry of it.

\section{Curvature tensor field of a statistical manifold}
\quad A statistical manifold $(M, g, \nabla)$  is defined as a manifold with a Riemannian metric $g$, a torsion-free affine connection $\nabla$ and a symmetric cubic form $Q(X, Y, Z)=(\nabla_{X}g)(Y, Z)$.
A new curvature equation which holds from nature characteristics of a statistical manifold is given and some properties associated with it are considered.\\
Let $M$ be a differential manifold and $g$ be a Riemannian metric on $M$. Let $\nabla$ be a torsion-free affine connection on $M$ and $\stackrel{*}\nabla$ be a dual connection of $\nabla$. Let $Q(X, Y, Z)=(\nabla_{X}g)(Y, Z)$ be a cubic form. 
\begin{Thot}\label{theo2_1}
 Let $(M, g, \nabla)$ be a statistical manifold. Then the curvature tensor field $R$ of $\nabla$ satisfies that
 \begin{equation}
 g(R(X, Y)Z, W)+g(R(Y, X)W, Z)=g(R(Z, W)X, Y)+g(R(W, Z)Y, X)\tag{2.1}
 \end{equation}
 for all $X, Y, Z, W\in T(M)$. 
\end{Thot}
\textbf{Proof}. From Ricci equation, we have
\begin{equation*}
 (\nabla g)(Z, W; Y; X)-(\nabla g)(Z, W; X; Y)=-g(R(X, Y)Z, W)-g(R(X, Y)W, Z)
 \end{equation*}
 for all $X, Y, Z, W\in T(M)$. Since $Q(X, Y, Z)=(\nabla g)(Y, Z; X)$, we have
\begin{equation*}
 g(R(X,Y)Z,W)+g(R(X,Y)W,Z)=(\nabla Q)(X,Z,W;Y)-(\nabla Q)(Y,Z,W;X)
\end{equation*}
On the other hand, from the definition of curvature tensor field and Bianchi equation , we have
\begin{eqnarray}
 && g(R(X,Y)Z,W)=-g(R(Y,X)Z,W)\nonumber\\
 && g(R(X,Y)Z,W)+g(R(Y,Z)X,W)+g(R(Z,X)Y,W)=0 \nonumber
\end{eqnarray}
for all $X, Y, Z, W\in T(M)$.\\
Hence we have 
\begin{eqnarray}
&& g(R(X,Y)Z,W)=-g(R(Y,Z)X,W)-g(R(Z,X)Y,W)\nonumber \\ 
&& =-g(R(Y,Z)W,X)+(\nabla Q)(Z,W,X;Y)-(\nabla Q)(Y,W,X;Z)\nonumber \\ 
&& \quad +g(R(Z,X)W,Y)+(\nabla Q)(X,W,Y;Z)-(\nabla Q)(Z,W,Y;X)\nonumber \\
\nonumber
\end{eqnarray}
\begin{eqnarray}
&& =-g(R(Z,W)Y,X)-g(R(W,Y)Z,X)+(\nabla Q)(Z,W,X;Y)-(\nabla Q)(Y,W,X;Z)\nonumber \\
&& \quad -g(R(X,W)Z,Y)-g(R(W,Z)X,Y)+(\nabla Q)(X,W,X;Z)-(\nabla Q)(Z,W,Y;X)\nonumber \\
&& =g(R(Z,W)X,Y)+(\nabla Q)(X,W,Y;Z)-(\nabla Q)(Z,X,Y;W)\nonumber \\ 
&& \quad +g(R(W,Y)X,Z)+(\nabla Q)(Y,X,Z;W)-(\nabla Q)(W,X,Z;Y)\nonumber \\
&& \quad +(\nabla Q)(Z,W,X;Y)-(\nabla Q)(Y,W,X;Z)\nonumber \\
&& \quad +g(R(X,W)Y,Z)+(\nabla Q)(W,Y,Z;X)-(\nabla Q)(X,Y,Z;W)\nonumber \\
&& \quad +g(R(Z,W)X,Y)+(\nabla Q)(X,W,Y;Z)-(\nabla Q)(Z,W,Y;X)\nonumber \\
&& =2g(R(Z,W)X,Y)-g(R(X,Y)Z,W)\nonumber \\
&& \quad +[(\nabla Q)(W,Y,Z;X)-(\nabla Q)(Z,W,Y;X)-(\nabla Q)(Y,Z,W;X)]\nonumber \\
&& \quad +[(\nabla Q)(X,Z,W;Y)+(\nabla Q)(Z,W,X;Y)-(\nabla Q)(W,X,Z;Y)]\nonumber \\
&& \quad +[(\nabla Q)(W,X,Y;Z)+(\nabla Q)(X,W,Y;Z)-(\nabla Q)(Y,W,X;Z)]\nonumber \\
&& \quad +[(\nabla Q)(Y,X,Z;W)-(\nabla Q)(Z,X,Y;W)-(\nabla Q)(X,Y,Z;W)]\nonumber 
\end{eqnarray}
for all $X, Y, Z, W\in T(M)$.\\
From the symmetry of $Q(X,Y,Z)$, we have
\begin{eqnarray}
&& g(R(X,Y)Z,W)=2g(R(Z,W)X,Y)-g(R(X,Y)Z,W)\nonumber \\ 
&& \quad -(\nabla Q)(Y,Z,W;X)-(\nabla Q)(Z,W,X;Y)+(\nabla Q)(W,X,Y;Z)-(\nabla Q)(X,Y,Z;W)\nonumber \\ 
&& =2g(R(Z,W)X,Y)-g(R(X,Y)Z,W)\nonumber \\ 
&& \quad +g(R(X,Y)Z,W)-g(R(Y,X)W,Z)-g(R(Z,W)X,Y)+g(R(W,Z)Y,X)\nonumber \\
&& =g(R(Z,W)X,Y)-g(R(Y,X)W,Z)+g(R(W,Z)X,Y)\nonumber 
\end{eqnarray}
for all $X, Y, Z, W\in T(M)$.\\
Arranging last equation, we obtain (2.1).$\quad\Box$ \\
In a statistical manifold $(M, g, \nabla)$, a torsion-free affine connection $\stackrel{*}\nabla$, called as a dual connection of $\nabla$, is naturally defined and dual structure $(g,\nabla,\stackrel{*}\nabla)$ is given.
Hence we can also denote a statistical manifold $(M,g,\nabla,\stackrel{*}\nabla)$ .
When $R$ and $\stackrel{*}{R}$ are curvature tensor field of $\nabla$ and dual connection $\stackrel{*}\nabla$ of $\nabla$, respectively, the following fact is satisfied.
\begin{Prop}\label{pro2_1} $([1])$
  In a statistical manifold $(M,g,\nabla,\stackrel{*}\nabla)$, we have 
 \begin{equation}
 g(R(X, Y)Z, W)+g(\stackrel{*}{R}(X, Y)W, Z)=0\quad \forall X, Y, Z, W\in T(M) \tag{2.2}
 \end{equation}
\end{Prop}
By using proposition 2.1, theorem 2.1 can be also expressed as follows:
\begin{align}
&& g(R(X, Y)Z, W)+g(\stackrel{*}{R}(X, Y)Z, W)=g(R(Z, W)X, Y)+g(\stackrel{*}{R}(Z, W)X, Y)\tag{2.3}\\
&& \forall X, Y, Z, W\in T(M) \nonumber
\end{align}
,where $R$ and $\stackrel{*}{R}$ are curvature tensor fields of $\nabla$ and dual connection $\stackrel{*}\nabla$ of $\nabla$, respectively.
\begin{Coro}\label{coro2_1}
  Let  $(M,g,\nabla,\stackrel{*}\nabla)$ be a statistical manifold and $Ric$ and $\stackrel{*}{Ric}$ be Ricci curvature tensor fields of $\nabla$ and dual connection $\stackrel{*}\nabla$ of $\nabla$, respectively. Then we have 
 \begin{equation}
 Ric(Y, Z)+\stackrel{*}{Ric}(Y, Z)=Ric(Z, Y)+\stackrel{*}{Ric}(Z, Y)\tag{2.4}
 \end{equation}
 for all $Y, Z\in T(M)$.
\end{Coro}
\textbf{Proof}. Since $Ric(Y,Z)=tr_{g}\{(X,W)\mapsto g(R(X,Y)Z,W)\}$ holds, contracting (2.3), we have 
\begin{eqnarray}
&& Ric(Y, Z)+\stackrel{*}{Ric}(Y, Z)\nonumber \\ 
&& \quad =tr_{g}\{(X,W)\mapsto g(R(X,Y)Z,W)\}+tr_{g}\{(X,W)\mapsto g(\stackrel{*}{R}(X,Y)Z,W)\}\nonumber \\ 
&& \quad =tr_{g}\{(X,W)\mapsto g(R(Z,W)X,Y)\}+tr_{g}\{(X,W)\mapsto g(\stackrel{*}{R}(Z,W)X,Y)\nonumber 
\end{eqnarray}
On the other hand, since $R(Z,W)X=-R(W,Z)X$ and 
\begin{equation*}
g(R(Z, W)X, Y)=-g(\stackrel{*}{R}(W, Z)X, Y)
\end{equation*}
hold, we have
\begin{eqnarray}
&& tr_{g}\{(X,W)\mapsto g(R(Z, W)X, Y)\}+tr_{g}\{(X,W)\mapsto g(\stackrel{*}{R}(Z, W)X, Y)\}\nonumber \\ 
&& \quad =tr_{g}\{(X,W)\mapsto g(\stackrel{*}{R}(W,Z)Y,X)\}+tr_{g}\{(X,W)\mapsto g(R(W,Z)Y,X)\}\nonumber\\
&& \quad =\stackrel{*}{Ric}(Z, Y)+Ric(Z, Y)\nonumber
\end{eqnarray}
Hence we obtain (2.4).$\quad\Box$ 
\begin{Def}\label{def2_1} $([7])$
Let $(M,g,\nabla)$  be a statistical manifold and $R$ be a curvature tensor field of $\nabla$. $(M,g,\nabla)$ is called a statistical manifold of constant curvature, if there is a constant $K$ satisfying that 
\begin{equation*}
R(X, Y)Z=K \{ g(Y,Z)X-g(X,Z)Y \},\qquad \forall X, Y, Z\in T(M)
\end{equation*}
\end{Def} 
The following facts hold.
\begin{Prop}\label{pro2_2}
 A statistical manifold of constant curvature $(M,g,\nabla)$ is a conjugate symmetric statistical manifold, that is , we have 
 \begin{equation}
 R(X, Y)Z=\stackrel{*}{R}(X, Y)Z, \qquad\forall X, Y, Z\in T(M) \tag{2.5}
 \end{equation}
where $R$ and $\stackrel{*}{R}$ are curvature tensor fields of $\nabla$ and dual connection $\stackrel{*}\nabla$ of $\nabla$, respectively.
\end{Prop}
\begin{Coro}\label{coro2_2}
 Let $(M,g,\nabla)$ be a statistical manifold of constant curvature. Then we have 
 \begin{equation*}
 Ric(Y, Z)=\stackrel{*}{Ric}(Y, Z),\qquad \forall Y, Z\in T(M)
 \end{equation*}
,where $Ric$ and $\stackrel{*}{Ric}$ are Ricci curvature tensor fields of $\nabla$ and dual connection $\stackrel{*}\nabla$ of $\nabla$, respectively.
\end{Coro}
\section{Conjugate Ricci-symmetry of a statistical manifold}
\quad In this section, we define conjugate Ricci-symmetry of a statistical manifold and consider a relation between equiaffine structure of a statistical manifold and it. We also give a relation between a statistical manifold of constant curvature and equiaffine structure of a statistical manifold.
\begin{Prop}\label{pro3_1}
 A statistical manifold $(M,g,\nabla)$ has an equiaffine structure if and only if a dual statistical manifold $(M,g,\stackrel{*}\nabla)$ has an equiaffine structure. 
\end{Prop}
Since a statistical manifold of constant curvature is a conjugate symmetric statistical manifold from Proposition 2.2, we have the following:
\begin{Prop}\label{pro3_2}
 If $(M,g,\nabla)$ is a statistical manifold of constant curvature, $\nabla^{(\alpha)}$ is equiaffine for all $\alpha\in\mathbf{R}$.
\end{Prop}
We consider a sufficient condition, which is weaker than conjugate symmetry , for $\nabla$ to be equiaffine.
\begin{Def}\label{def3_1}
Let $(M,g,\nabla,\stackrel{*}\nabla)$ be a statistical manifold and $Ric$ and $\stackrel{*}{Ric}$ are Ricci curvature tensor fields of $\nabla$ and dual connection $\stackrel{*}\nabla$ of $\nabla$, respectively. We say that a statistical manifold $(M,g,\nabla,\stackrel{*}\nabla)$ is conjugate Ricci-symmetric, if we have  
 \begin{equation}
 Ric(Y, Z)=\stackrel{*}{Ric}(Y, Z)\tag{3.1}
 \end{equation}
for all $Y, Z\in T(M)$.
\end{Def} 
\begin{Thot}\label{theo3_1}
 If a statistical manifold $(M,g,\nabla,\stackrel{*}\nabla)$ is conjugate Ricci-symmetric , $\nabla$ is equiaffine.
\end{Thot}
\textbf{Proof}. From corollary 2.1, we have
\begin{equation*}
 Ric(Y, Z)+\stackrel{*}{Ric}(Y, Z)=Ric(Z, Y)+\stackrel{*}{Ric}(Z, Y)
\end{equation*}
 for all $Y, Z\in T(M)$ and since from the condition of theorem we have
\begin{equation*}
 Ric(Y, Z)=\stackrel{*}{Ric}(Y, Z),\quad Ric(Z, Y)=\stackrel{*}{Ric}(Z, Y) \quad\forall Y, Z\in T(M)
\end{equation*}
Hence $\nabla$ is equiaffine.$\quad\Box$ 
\begin{Coro}\label{coro3_1}
 If a statistical manifold $(M,g,\nabla,\stackrel{*}\nabla)$ is conjugate Ricci-symmetric, $\nabla^{(\alpha)}$ is equiaffine for all $\alpha\in\mathbf{R}$.
\end{Coro}
\begin{lemm}\label{lemm3_1}
If there exist an $\alpha_0\in\mathbf{R}(\alpha_{0}\neq 0)$ such that $Ric^{(\alpha_{0})}=Ric^{(-\alpha_{0})}$, we have $Ric^{(\alpha)}=Ric^{(-\alpha)}$ for all $\alpha\in\mathbf{R}$, where $Ric^{(\alpha)}$ is a Ricci curvature tensor of $\alpha $-connection of $\nabla^{(\alpha)}$.
\end{lemm} 
\textbf{Proof}. We have
\begin{equation}
Ric^{(\alpha)}(Y, Z)
=\frac{1+\alpha}{2}Ric(Y,Z)+\frac{1-\alpha}{2}\stackrel{*}{Ric}(Y, Z)+\frac{1-\alpha^2}{4}(trK(K(\cdot ,Z),Y)-trK(\cdot ,(K(Y,Z)))\tag{3.2}
\end{equation}
for all $Y, Z\in T(M)$, where $K(Y,Z)=\stackrel{*}\nabla_{Y} Z-\nabla_{Y} Z$ is the ''difference tensor''.
From (3.2), we have
\begin{equation}
Ric^{(\alpha)}(Y,Z)-Ric^{(-\alpha)}(Y,Z)=\alpha(Ric(Y,Z)-\stackrel{*}{Ric}(Y,Z))
\tag{3.3}
\end{equation}
(3.3) implies that if there exist an $\alpha_0\in\mathbf{R}(\alpha_{0}\neq 0)$ such that $Ric^{(\alpha_{0})}=Ric^{(-\alpha_{0})}$, we have $Ric(Y,Z)=\stackrel{*}{Ric}(Y,Z)$. (3.3) also implies that if $Ric(Y,Z)=\stackrel{*}{Ric}(Y,Z)$ holds, we have $Ric^{(\alpha)}=Ric^{(-\alpha)}$ for all $\alpha\in\mathbf{R}(\alpha\neq 0)$. On the other hand if $\alpha=0$ holds, we have $Ric^{(\alpha)}=Ric^{(-\alpha)}=\stackrel{\circ }{Ric}$.$\quad\Box$\\ 
Lemma 3.1 and corollary 3.1 imply the following fact.
\begin{Thot}\label{theo3_2}
If there exist an $\alpha_0\in\mathbf{R}(\alpha_{0}\neq 0)$ such that $Ric^{(\alpha_{0})}=Ric^{(-\alpha_{0})}$, $\nabla ^{(\alpha)}$ is equiaffine for all $\alpha\in\mathbf{R}$.
\end{Thot}
Then direct calculation shows:
\begin{Prop}\label{pro3_3}
 A 2-dimensional statistical manifold $(M,g,\nabla,\stackrel{*}\nabla)$ is conjugate symmetric if and only if it is conjugate Ricci-symmetric.
\end{Prop}
\section{Condition for conjugate symmetry and conjugate Ricci-symmetry of a statistical manifold to coincide}
\quad In this section, we define a statistical manifold $(M,g,Q)$  with a recurrent metric  and consider a condition for conjugate symmetry and conjugate Ricci-symmetry of a statistical manifold to coincide there.
\begin{Def}\label{def4_1}
Let $(M,g,\nabla,\stackrel{*}\nabla)$ be a statistical manifold. $(M,g,\nabla,\stackrel{*}\nabla)$ is called a statistical manifold with a recurrent metric, if there exist a 1-form $\omega $ on $M$ satisfying that 
\begin{equation}
(\nabla g)(Y, Z;X)=\omega(X)g(Y,Z)+\omega(Y)g(Z,X)+\omega(Z)g(X,Y)\tag{4.1} 
\end{equation}
for all $X, Y, Z\in T(M)$.
\end{Def} 
\begin{Def}\label{def4_2}
Let $(M,g,\nabla,\stackrel{*}\nabla)$ be a statistical manifold and $(M,h,\stackrel{\circ }\nabla)$ be a Riemannian manifold. We say that $(M,g,\nabla,\stackrel{*}\nabla)$ and $(M,h,\stackrel{\circ }\nabla)$ are $\alpha$-conformal, if there exist a $\varphi \in C^{\infty }(M)$ satisfying that 
\begin{eqnarray}
&& 1^{\circ } \quad g(X,Y)=e^{\varphi }h(X,Y)\nonumber\\
&& 2^{\circ } \quad \nabla_{X}Y=\stackrel{\circ }\nabla_{X}Y+\frac{1-\alpha}{2}(d\varphi(X)Y+d\varphi(Y)X)-\frac{1+\alpha}{2}h(X,Y)grad_{h}\varphi\nonumber
\end{eqnarray}
for all $X, Y\in T(M)$.
\end{Def} 
\begin{lemm}\label{lemm4_1}
A statistical manifold $(M,g,\nabla,\stackrel{*}\nabla)$, which is $\alpha$-conformal with a Riemannian manifold  $(M,h,\stackrel{\circ }\nabla)$ , is a statistical manifold with a recurrent metric.
\end{lemm} 
\begin{Coro}\label{coro4_1}
 A statistical manifold $(M,g,\nabla)$ is $\alpha$-conformal with a Riemannian manifold  $(M,h,\stackrel{\circ }\nabla)$ if and only if a dual statistical manifold $(M,g,\stackrel{*}\nabla)$ is $(-\alpha)$-conformal with a Riemannian manifold $(M,h,\stackrel{\circ }\nabla)$.
\end{Coro}
A conjugate symmetry and a conjugate Ricci-symmetry of a statistical manifold $(M,g,Q)$ with a recurrent metric coincide under some conditions.
\begin{Thot}\label{theo4_1}
Let $(M,g,Q)$ be a $n$-dimensional statistical manifold with a recurrent metric. Let 
\begin{equation*}
(\nabla g)(Y, Z;X)=\omega(X)g(Y,Z)+\omega(Y)g(Z,X)+\omega(Z)g(X,Y)
\end{equation*}
be a cubic form , where $\omega$ is a closed form. Then $(M,g,Q)$ is conjugate symmetry if and only if it is conjugate Ricci-symmetry.
\end{Thot}
\textbf{Proof}. We have
\begin{equation*}
\stackrel{*}{Ric}(Y,Z)-Ric(Y,Z)=tr_{g}{(X,W)\mapsto [(\stackrel{\circ }\nabla_{X} Q)(Y,Z,W)-(\stackrel{\circ}\nabla_{Y} Q)(Y,Z,W)]}
\end{equation*}
for all $X, Y, Z\in T(M)$.
On the other hand, since we have
\begin{equation*}
(\nabla g)(Y, Z;X)=\omega(X)g(Y,Z)+\omega(Y)g(Z,X)+\omega(Z)g(X,Y)
\end{equation*}
and $\omega$ is closed form, that is, we have 
\begin{equation*}
X(\omega(Y))=Y(\omega(X)),\quad (\stackrel{\circ}\nabla_{X} \omega)(Y)=(\stackrel{\circ}\nabla_{Y} \omega)(X)
\end{equation*}
for all $X, Y\in T(M)$,we have
\begin{equation*}
(\stackrel{\circ}\nabla_{X} Q)(Y,Z,W)-(\stackrel{\circ}\nabla_{Y} Q)(X,Z,W)
\end{equation*}
\begin{equation}
 =g(W,Y)(\stackrel{\circ}\nabla_{X} \omega)(Z)-g(W,X)(\stackrel{\circ}\nabla_{Y} \omega)(Z)+g(Y,Z)(\stackrel{\circ}\nabla_{X} \omega)(W)-g(X,Z)(\stackrel{\circ}\nabla_{Y} \omega)(W)\tag{4.2}
\end{equation}
and
\begin{eqnarray}
&& \stackrel{*}{Ric}(Y,Z)-Ric(Y,Z)=tr_{g}{(X,W)\mapsto [(\stackrel{\circ }\nabla_{X} Q)(Y,Z,W)-(\stackrel{\circ}\nabla_{Y} Q)(Y,Z,W)]}\nonumber\\
&& \quad =g(Y,Z)(\stackrel{\circ }\nabla \omega )(\cdot ,\cdot )-n(\stackrel{\circ }\nabla_{Y} \omega )(Z)\nonumber
\end{eqnarray}
Hence, if $(M,g,Q)$ is conjugate Ricci-symmetry, we have
\begin{equation*}
(\stackrel{\circ }\nabla_{Y} \omega )(Z)=\frac{1}{n} g(Y,Z)(\stackrel{\circ }\nabla \omega )(\cdot ,\cdot )
\end{equation*}
for all $Y, Z\in T(M)$ and, substituting this equation to (4.2), we have
\begin{equation*}
(\stackrel{\circ}\nabla_{X} Q)(Y,Z,W)-(\stackrel{\circ}\nabla_{Y} Q)(X,Z,W)=0
\end{equation*}
for all $X,Y, Z,W\in T(M)$,that is, $(M,g,Q)$ is conjugate symmetry.$\quad\Box$\\
Theorem 4.1 and lemma 4.1 imply the following fact.
\begin{Coro}\label{coro4_2}
 A statistical manifold $(M,g,\nabla,\stackrel{*}\nabla)$, which is $\alpha$-conformal with a Riemannian manifold $(M,h,\stackrel{\circ }\nabla)$, is conjugate symmetry if and only if it is conjugate Ricci-symmetry.
\end{Coro}

\end{document}